\documentclass[preprint, 12pt]{article}
 \usepackage{latexsym}
 \usepackage{graphicx}
 \usepackage{amssymb}
 \title{Edge Decompositions of Hypercubes by Paths}
\begin{document}
% \begin{frontmatter}
 
 \author{David Anick\\
 \small Laboratory for Water and Surface Studies\\[-.08ex]
 \small Tufts University \\[-0.8ex]
 \small Department of Chemistry \\[-0.8ex]
 \small 62 Pearson Rd. \\[-0.8ex]
 \small Medford, MA 02155\\[-0.8ex]
 \small \texttt{david.anick@rcn.com}\\
 \small and\\
 Mark Ramras$^*$\\
 \small Department of Mathematics\\[-0.8ex]
 \small Northeastern University\\[-0.8ex]
  \small Boston, MA 02115, USA\\[-0.8ex]
 \small \texttt{m.ramras@neu.edu}\\
\small (617)-373-5651\\[-0.8ex]
 \small $^*$Corresponding author}
 \maketitle

 \begin{abstract}

 Many authors have investigated edge decompositions of graphs by the edge sets of isomorphic copies of special subgraphs.  For $q$-dimensional hypercubes $Q_q$  various researchers have done this for certain trees, paths, and cycles.  In this paper we shall say that ``$H$ divides $G$"  if $E(G)$ is the disjoint union of $\{E(H_i)\,|\,H_i \simeq H\}$.   Our main result is that for $q$ odd and $q<2^{32}$, the path of length $m, P_m,$ divides $Q_q$ if and only if $m\leq q$ and $m\,|\,q\cdot 2^{q-1}$. 

 \end{abstract}
%\begin{keyword} 
 \newpage
\section{Introduction}
   
   Edge decompositions of graphs by subgraphs has a long history.  For example, there is a Steiner triple system of order $n$ if and only if the complete graph $K_n$ has an edge-decomposition 
 by $K_3$.  In 1847 Kirkman \cite{Kirk1} proved that for a Steiner triple system to exist it is necessary that $n\equiv 1 \pmod {6}$ or $n\equiv 3 \pmod {6}$.  In 1850 he \cite{Kirk2} proved the converse holds also.
 
 \newtheorem{theorem}{Theorem}
 \begin{theorem}\label{Steiner}   
 A Steiner system of order $n\geq 3$ exists if and only if $n\equiv 1 \pmod {6}$ or $n\equiv 3 \pmod {6}$.
 \end{theorem}
 
 In more modern times (1964) G. Ringel \cite{Ring2} stated the following conjecture, which is still open.\\
 \vspace{.1in}

\noindent  {\bf Ringel's Conjecture}\\
 If $T$ is a fixed tree with $m$ edges then $K_{2m+1}$ is edge-decomposable into $2m+1$ copies of $T$.
 \vspace{.1in}

Still more recently, the $n$-dimensional hypercube graph $Q_n$ has been studied extensively, largely because of its usefulness as the architecture for distributed parallel processing supercomputers \cite{Lei}.  Communication problems such as ``broadcasting" in these networks (see \cite{Jo-Ho}, \cite{Bruck}) have led to research on constructions of maximum size families of edge-disjoint spanning trees (maximum is $\lfloor{n/2}\rfloor$ for $Q_n$ \cite{Barden et al}; see \cite{Ku-Wang-Hung} for results on more general product networks.)   Fink \cite{Fink} and independently Ramras \cite{Ram1} proved that $Q_n$ could be decomposed into $2^{n-1}$ isomorphic copies of any tree on $n$ edges.  Wagner and Wild \cite{Wag-Wild} proved that $Q_n$ is edge-decomposable into $n$ copies of a tree on $2^{n-1}$ edges.  Horak, Siran, and Wallis \cite{Hor-Sir-Wal} showed that $Q_n$ has an edge decomposition by isomorphic copies of any graph G with $n$ edges each of whose blocks is either an even cycle or an edge. Ramras \cite{Ram2} proved that for a certain class of trees on $2n$ edges, isomorphic copies of these trees edge-decompose $Q_n$.    Other researchers have demonstrated edge decompositions by Hamiltonian cycles for Cartesian products of cycles \cite{Ring1}, \cite{AS}, \cite{E-ZVE}.  Song \cite{Song} applies a different construction of this to even-dimensional hypercubes.

We concentrate in this work on the important question of edge decompositions of hypercubes into paths of equal length.  Literature on this specific question is not extensive.  The cases of $n$ odd and $n$ even are very different, with the theory of edge decompositions of $Q_n$ for $n$ even being dominated by Hamiltonian cycle considerations as noted above.  Mollard and Ramras \cite{Mol-Ram} found edge decompositions of $Q_n$ into copies of $P_4$, the path on 4 edges, for all $n\geq 5$.  Our principal result goes far beyond that: for $n < 2^{32}$ we answer the general question of when $Q_n$ for $n$ odd can be edge decomposed into length-$m$ paths.  The method of proof involves construction of two new graph-theoretic operations that may have wide applicability to edge decomposition studies.

\section{Notation and Preliminaries}
\newtheorem{definition}{Definition}
\begin{definition}\label{<D}
For graphs $H$ and $G$ we say that $H$ \underline{divides} $G$ if there is a collection of subgraphs $\{H_i\}$ each 
isomorphic to $H$ ($H_i \simeq H$ for all $i$) for which $E(G)$ is the disjoint union of $\{E(H_i)\}$.
\end{definition}

\noindent {\bf Notation}  We shall denote ``$H$ divides $G$" by $H <_D G$ since the relation  $<_D$ is clearly reflexive and transitive and thus a partial order.

For the $q$-dimensional hypercube $Q_q$ the vertices are the $2^q\, q$-tuples of $0$'s and $1$'s.  $V(Q_q)$ has an additive structure of $\mathbb{Z}_2^q$.  The edge set $E(Q_q)$ consists of those (unordered) pairs of vertices that differ in exactly one 
coordinate.  The group $\mathbb{Z}_2^q$ acts on the set of edges in the obvious way: for  $\gamma\in\mathbb{Z}_2^q$ and $e=\{\alpha,\alpha^{\,\prime}\}$ an unordered pair representing an edge of $E(Q_q)$, $\gamma+e$ will denote the edge 
$\{\gamma+\alpha,\gamma+\alpha^{\,\prime}\}$.

The parity of a q-tuple $\alpha=(a_1, \ldots , a_q)\in\mathbb{Z}_2^q$ is 
$\rho(\alpha)=a_1+\ldots , a_q,$ defined$\pmod 2.$  Let $B_q$ be the subgroup of $V(Q_q)$ consisting of those $q$-tuples with parity $0$.  For $q\ge 1$ clearly $|\,B_q\,|=|\,V(Q_q)\,|/2=2^{q-1}.$
 
 Given an integer $j, 1\leq j\leq q$, and a vertex $\alpha=(a_1, \ldots , a_q)\in V(Q_q)$, some helpful notations are as follows. Let 
 $$\overline{j}\cdot \alpha=(a_1, \ldots , 1+a_j, \ldots , a_q)$$
   {\it i.e.} alter $a_j$ only.  Let
 $$j^0\cdot \alpha=(a_1, \ldots , c, \ldots , a_q),$$
 where $c=\rho(\alpha)+a_j$.  The idea of $j^0$ is 
``alter the $j^{th}$ coordinate if necessary so that the parity is $0$".  It should be obvious that $j^0\cdot \alpha=j^0\cdot (\overline{j}\cdot \alpha)\in B_q.$  Likewise put
$$j^1\cdot \alpha=\overline{j}\cdot (j^0\cdot \alpha),$$ 
{\it i.e.}  alter the $j^{th}$ coordinate if necessary so that the parity is 1.   
Notice that $\{\alpha, \overline{j}\cdot \alpha\}$ is an edge of $Q_q$ and that 
$\{j^0\cdot \alpha, j^1\cdot \alpha\}$ is the same edge.  
Our notation for this edge is $\hat{j}\cdot \alpha.$  
Then $\hat{j}$ is compatible with the $\mathbb{Z}_2^q$-action, {\it i.e.} 
$\hat{j}\cdot(\gamma+\alpha)=\gamma+\hat{j}\cdot\alpha.$  Clearly 
$\hat{j}\cdot \alpha=\hat{j}\cdot (j^0\cdot\alpha)=\hat{j}\cdot (j^1\cdot\alpha)
=\hat{j}\cdot (\overline{j}\cdot\alpha).$

The path $P_q$ of length $q$ is a graph with with vertex set $\{0, 1, \ldots , q\}$ and edge set $\{\hat{1},\ldots,\hat{q}\}$, $\hat{k}$ denoting the edge joining $k-1$ and $k$.  
We define graph embeddings $f_{\gamma}: P_q\longrightarrow Q_q,$ for $\gamma\in B_q,$
as follows.  For $0 \leq k \leq q$ let $$1^k0^{q-k}=(\underbrace{1, \ldots , 1}_{k \,1{\rm 's}}, \underbrace{0, \ldots , 0}_{q-k\, 0{\rm 's}})\in V(Q_q)$$
and set
$$f_{\gamma}(k) = 1^k0^{q-k}+\gamma.$$
Notice that in $E(Q_q)$,
$$f_{\gamma}(\hat{k}) =\hat{k}\cdot(1^k0^{q-k}+\gamma)=1^k0^{q-k}+\hat{k}\cdot\gamma
=1^{k-1}0^{q-k+1}+\hat{k}\cdot\gamma.$$

The family $\{f_{\gamma}\}$ provides $|\,B_q\,|=2^{q-1}$ ways of embedding $P_q$ in $Q_q$, and $P_q$ has $q$ edges, so altogether the $\{f_{\gamma}\}$ send $q\cdot 2^{q-1}$ 
edges to $Q_q$ while $|\,E(Q_q)\,|=q\cdot 2^{q-1}.$  Therefore if the $\{f_{\gamma}\}$ cover
$E(Q_q)$ then they cover each edge just once, {\it i.e.} the path images of $\{f_{\gamma}\}$
are pairwise edge-disjoint.  To see that this is the case, let $e=\{\alpha,\alpha^{\,\prime}\}$ denote any edge of $Q_q;$ then $e=\hat{k}\cdot\alpha$ where the unique coordinate that
differs between $\alpha$ and $\alpha^{\,\prime}$ is the $k^{\textrm{th}}.$  Put
$\gamma=k^0\cdot(\alpha+1^k0^{q-k}).$ and observe that 
$f_{\gamma}(\hat{k})=\hat{k}\cdot\gamma=e.$  We have proved

\newtheorem{lemma}{Lemma}
\begin{lemma}\label{fgammabijects}
The family of graph embeddings $\{f_{\gamma}: P_q\longrightarrow Q_q\,|\,\gamma\in B_q\}$ defines a partition of $E(Q_q)$ into edge-disjoint paths indexed by $B_q.$
\end{lemma}

(As mentioned in the Introduction, a more general result, for {\em all} trees on $q$ edges, appears in \cite{Fink} and in \cite{Ram1}.)

The results in the next lemma are also in \cite{Ram1} but we include short proofs here so this article can be self-contained.

\begin{lemma}\label{P2 divides Q3}
{\rm (a)}  $P_2 <_D Q_3$.\\
%%{\rm (b)}  $P_4$ divides $Q_5$.\\
{\rm (b)}  If $P_{2^m} <_D Q_q$, where $q$ is odd, then $q\geq 2^m$.  
\end{lemma}
\vspace{.1in}

\noindent {\em Proof.}  (a)  $Q_3$ may be viewed as an inner $Q_2$ joined to an outer $Q_2$ via a perfect matching.  Decompose the inner $Q_2$ into 2 edge-disjoint $P_2$'s.  Each of the remaining 8 edges decompose into 4 $P_2$'s, with one edge of the outer $Q_2$ joined to an incident matching edge.
\vspace{.1in}

 (b)  Every vertex of $Q_q$ has odd degree, so at every vertex at least one embedded path must start or end there.  So there must be at least $|\,V(Q_q)\,|/2$ paths, {\it i.e.} $q\cdot 2^{q-1}/2^m\geq 2^q/2$, which implies that  $q\geq 2^m.$  \hfill  $\Box$

\section{Stretched Graphs}

\begin{definition}\label{Definition*}
Let $G$ be a graph and let $m$ be a positive integer.  The \underline{m-stretch} of $G,$ 
denoted $m*G,$ is the graph obtained by replacing each edge of $G$ by a path of length $m$, 
so that these paths are internally vertex-disjoint.
\end{definition} 

\begin{lemma}\label{stretchprops} 
 {\rm (a)}  $1*G \simeq G$ for any graph $G$.\\
 {\rm (b)}  $|\,E(m*G)\,|=m|\,E(G)\,|$.\\
 {\rm (c)} $ |\,V(m*G)\,|=(m-1)|\,E(G)\,|+|\,V(G)\,|$.\\
 {\rm (d)}  $m_1*(m_2*G)\simeq (m_1m_2)*G$.\\
 {\rm (e)}  If $H <_D G$, then $m*H <_D m*G$.\\
 {\rm (f)}  $m*P_q\simeq P_{mq}$.
 
 \end{lemma}
 The proofs are trivial.
 
The importance for hypercubes of stretched graphs comes from

\begin{theorem}\label{m*Qq divides Qmq}
$m*Q_q <_D Q_{mq}$ for any $m \geq 1, q \geq 1.$ 
 \end{theorem}

For example, from this and Lemmas \ref{P2 divides Q3}(a) and \ref{stretchprops}(e,f) it follows easily that $P_6=3*P_2 <_D 3*Q_3$, which divides $Q_9$.  By transitivity of divisibility, one obtains $P_6  <_D Q_9$, which is already a new result.  To prove Theorem \ref{m*Qq divides Qmq}, the cases of $m$ odd and $m=2$ are considered separately.  It should be clear from Lemma \ref{stretchprops}(d,e) that if $m_1*Q_q <_D Q_{m_1q}$ for any $q$ and if $m_2*Q_q <_D Q_{m_2q}$ for any $q,$ then $m_1m_2*Q_q <_D Q_{m_1m_2q}$ for any $q,$ so the cases of $m$ odd and $m=2$ suffice.
\vspace{.1in}
 
  \newtheorem{proposition}{Proposition}
 \begin{proposition}\label{case m odd}
 $m*Q_q <_D Q_{mq}$ for $m$ odd, $q\geq 1$.
 \end{proposition}
 
 \noindent  Before jumping into the proof, let us establish notations for vertices 
and edges of $Q_{mq}$ and $m*Q_q$.  
We consider a vertex of $Q_{mq}$ to consist of $q$ vectors of length $m$,
 (view $Q_{mq}$ as $Q_m^q$) {\it i.e.} $\alpha=(\alpha^{(1)}, \alpha^{(2)}, \ldots , \alpha^{(q)}), \alpha^{(k)}\in V(Q_m)=\mathbb{Z}_2^m$.  Like before $0^m$ is $(0,\ldots, 0)\in \mathbb{Z}_2^m$ and $1^m$ is $(1,\ldots,1)\in \mathbb{Z}_2^m$.
 
 Notations for $m*Q_q$ are as follows.  First, each vertex of $Q_q$ is carried over as a vertex into $m*Q_q$, so if $\alpha=(a_1, \ldots, a_q)\in V(Q_q)$, we also view $\alpha$ as a vertex of $m*Q_q$.  In addition, for each edge $\hat{j}\cdot \alpha \in E(Q_q)$, let $j^k\!:\!\alpha$ denote the $k^{\textrm{th}}$ point on the path that replaced $\hat{j}\cdot \alpha,$ where $0\leq k\leq m$.  We also identify $j^0\!:\!\alpha$ with $\alpha$, and $j^m\!:\!\alpha$ with $\overline{j}\cdot \alpha$ (which is the other endpoint of $\hat{j}\cdot \alpha)$.  Note that the edges of $m*Q_q$ connect $j^{k-1}\!:\!\alpha$ with $j^k\!:\!\alpha, k=1, \ldots, m$.  The points and edges can be counted coming from either end of the path, hence
 $$j^k\!:\!\alpha=j^{m+1-k}\!:\!(\overline{j}\cdot \alpha).$$
 So one must be careful that any definition involving $j^k\!:\!\alpha$ is independent of choice of notation.  One way to make the above notation unique for the vertices not inherited from $V(Q_q)$ is to apply it only to $\alpha\in B_q$.  Then
$$V(m*Q_q)=V(Q_q)\cup \{j^k\!:\!\alpha \,|\,1\leq k<m, 1\leq j \leq q, \alpha\in B_q\}.$$

\noindent {\em Proof} of Proposition \ref{case m odd}.\\
Let $\gamma=(\gamma^{(1)}, \gamma^{(2)}, \ldots, \gamma^{(q)})\in B_m^q\subseteq \mathbb{Z}_2^{mq},$ {\it i.e.} a vector where each length-$m$ subvector $\gamma^{(i)}$ has parity $0$.  Define embeddings $F_{\gamma}:m*Q_q\rightarrow Q_{mq}$
as follows.  If $\alpha=(a_1, \ldots, a_q)\in V(Q_q)$, put 
$$F_{\gamma}(\alpha)=(a_1^{m}, a_2^m, \ldots, a_q^m)+\gamma.$$
 If instead a vertex of $m*Q_q$ is $j^k\!:\!\alpha,$ with $\alpha\in B_q$, put
$$F_{\gamma}(j^k\!:\!\alpha)=(c^{(1)}, c^{(2)}, \ldots , c^{(q)})+\gamma,$$
 where $c^{(s)}=\left\{ \begin{array}{ll}
                                            a_s^m & {\rm if}\, s\neq j\\
                                            a_j^m+1^k0^{m-k} & {\rm for}\, s=j.
                                            \end{array} \right.$
Note that $F_{\gamma}(j^0\!:\!\alpha)=F_{\gamma}(\alpha)$ by this definition and likewise $F_{\gamma}(j^m\!:\!\alpha)=F_{\gamma}(\overline{j}\cdot \alpha)$, as needed for notational consistency and for $F_{\gamma}$ to send edges to edges.

We will show that the $\{F_{\gamma}\}$ comprise an edge partition of $Q_{mq}$ into copies of $m*Q_q$.  Now $|\,E(m*Q_q)\,|=m\cdot (q\cdot 2^{q-1})=mq\cdot 2^{q-1},$ and with $|\,B_m^q\,|=(2^{m-1})^q=2^{mq-q}$ embeddings, at most $(2^{mq-q})(mq\cdot 2^{q-1})=mq\cdot 2^{mq-1}$ edges will be covered by the union of their images.  But this is exactly $|\,E(Q_{mq})\,|$, so the $\{F_{\gamma}\}$ comprise an edge partition if and only if                                             
 $$\bigcup_{\gamma}F_{\gamma}(E(m*Q_q))\supseteq E(Q_{mq}),$$
 {\it i.e.} it suffices to show that every edge of $Q_{mq}$ is in the image of some $F_{\gamma}$.
 
 Let an edge of $Q_{mq}$ be written as $(\alpha^{(1)}, \ldots , \hat{k}\cdot \beta, \ldots , \alpha^{(q)})$, where $\alpha^{(s)}\in Q_m$ and $\beta\in B_m$.  The idea here is that the unique vertex that changes over the edge is at some position (call it $k$) of some length-$m$ segment (call it the $j^{\textrm{th}}$).  Put 
 $$\gamma^{(s)}=\rho(\alpha^{(s)})^m+\alpha^{(s)}\,\, {\rm for}\, s\neq j.$$
 Then $\gamma^{(s)}\in B_m$ because the parity of $m$ copies of either $0\, {\rm or}\, 1$ is (respectively) either $0$\, or\, $1$.  (Note:  This is the {\it only} place in the proof where the 
hypothesis is used that $m$ is odd.)  
 Put $$c=\rho(\alpha^{(1)})+\ldots , + \rho(\alpha^{(j-1)})+\rho(\alpha^{(j+1)})+ \ldots, +\rho(\alpha^{(q)})$$
and set
 $$\gamma^{(j)}=k^0\cdot (c^m+1^k0^{m-k}+\beta).$$
Putting $\gamma=(\gamma^{(1)}, \gamma^{(2)}, \dots , \gamma^{(q)}),$ 
we have $\gamma\in B_m^q.$  
Let $b_s=\left\{\begin{array}{lll} c & {\rm for}\,s=j & \\
                                                          \rho(\alpha^{(s)}) & {\rm for}\, s\neq j, &  \end{array} \right.$ and let $\beta=(b_1, \ldots , b_q).$  
Then $\beta\in B_q$ and 
$$F_{\gamma}(j^{k}\!:\!\beta)=(b_1^m+\gamma^{(1)}, \ldots, c^m+1^k0^{m-k}+\gamma^{(j)}, \ldots , b_q^m+\gamma^{(q)})$$
$$=(\alpha^{(1)}, \ldots , k^{c+k}\cdot \beta, \ldots , \alpha^{(q)}).$$
and likewise 
$$F_{\gamma}(j^{k-1}:\beta)=(b_1^m+\gamma^{(1)}, \ldots , c^m+1^{k-1}0^{m-k+1} +\gamma^{(j)}, \ldots , b_q^m+\gamma^{(q)})$$
$$=(\alpha^{(1)}, \ldots , k^{c+k-1}\cdot \beta, \ldots , \alpha^{(q)}).$$ 

This shows that the edge of $m*Q_q$ that links $j^k\!:\!\beta$ and $j^{k-1}\!:\!\beta$ is sent via $F_{\gamma}$ to the edge $(\alpha^{(1)}, \ldots , \hat{k}\cdot \beta, \ldots , \alpha^{(q)})$ of $Q_{mq}$.  This completes the proof of Proposition \ref{case m odd}.  \hfill  $\Box$

\begin{proposition}\label{2*Qq divides Q2q} 
$2*Q_q <_D Q_{2q}$, for any $q\geq 1$.
\end{proposition}

\noindent {\em Proof.}  Again, begin with notation for vertices and edges of $Q_{2q}$ and $2*Q_q$.  This time we set up $Q_{2q}$ slightly differently, namely we identify $Q_{2q}$ as $\mathbb{Z}_2^q\times \mathbb{Z}_2^q$ so $(\alpha, \beta)$ would be a typical vertex of $Q_{2q}$, where $\alpha , \beta \in V(Q_q)$.  The notation for $V(2*Q_q)$ is similar to before, but there is no need for the superscript `$k$' because $k=1$, and we will simply write 
$j\!:\!\alpha$ for the midpoint of $\hat{j}\cdot \alpha$.  Notice that $j\!:\!\alpha
=j\!:\!(\overline{j}\cdot \alpha)$.  A unique notation for $V(2*Q_q)$ is implicit in
$$V(2*Q_q)=V(Q_q)\cup \{j\!:\!\alpha\,|\,\alpha\in B_q, 1\leq j\leq q\}.$$
Let $j\#\alpha$ denote the unique edge of $2*Q_q$ connecting $\alpha$ and $j\!:\!(j^0\cdot \alpha)$.  For any $\gamma\in B_q$, define the embeddings
$F_{\gamma}^0, F_{\gamma}^1$ : $2*Q_q\rightarrow Q_{2q}$ by
$$F_{\gamma}^0(\alpha)=F_{\gamma}^1(\alpha)=(\alpha, \alpha+\gamma);$$
$$F_{\gamma}^0(j\!:\!\alpha)=(j^0\cdot \alpha, (j^1\cdot \alpha)+\gamma);$$
$$F_{\gamma}^1(j\!:\!\alpha)=(j^1\cdot \alpha, (j^0\cdot \alpha)+\gamma).$$
Note that $|\,E(2*Q_q)\,|=2|\,E(Q_q)\,|=2q\cdot 2^{q-1}=q\cdot 2^q$, and there are $2^{q-1}$ elements in $B_q$ and so $2\cdot 2^{q-1}=2^q$ embeddings.  Their combined images cover at most $(2^q)(q\cdot 2^q)=q\cdot 2^{2q}$ edges, and again 
$|\,E(Q_{2q})\,|=(2q)\cdot 2^{2q-1}=q\cdot 2^{2q}$ also, so the family
$\{F_{\gamma}^{\epsilon}\,|\, \epsilon=0 \textrm{ or } 1, \gamma\in B_q\,\}$ 
provides an edge partition of $Q_{2q}$ into copies of $2*Q_q$ if and only if
$$E(Q_{2q})\subseteq \bigcup_{\gamma, \epsilon}F_{\gamma}^{\epsilon}(E(2*Q_q)).$$

To prove this we break into two cases.  An edge $e$ of $Q_{2q}$ is either $(\hat{j}\cdot \alpha, \beta)$, where $\alpha\in B_q$ and $\beta\in \mathbb{Z}_2^q$, or $(\alpha, \hat{j}\cdot \beta),$ where $\alpha\in \mathbb{Z}_2^q, \beta\in B_q.$  Define $\gamma$ by $\gamma=j^0\cdot(\alpha+\beta)\in B_q$, and let $\tilde{\alpha}=\beta+\gamma\in B_q$.

Suppose the edge $e$ is 
$(\hat{j}\cdot \alpha, \beta)$, with $\alpha\in B_q$ and $\beta\in \mathbb{Z}_2^q.$
Divide further into two subcases.  If $\rho(\beta)=0$, then $\gamma=\alpha+\beta$ and $\tilde{\alpha}=\alpha\in B_q$.  
We have $$F_{\gamma}^1(j\!:\!\alpha)=(j^1\cdot \alpha, (j^0\cdot \alpha)+\gamma)=(\overline{j}\cdot\alpha, \beta)$$ while
$$F_{\gamma}^1(\alpha)=(\alpha, \alpha+\gamma)=(\alpha, \beta).$$
Hence $F_{\gamma}^1$ carries $j\#\alpha$ to $(\hat{j}\cdot \alpha, \beta).$

Otherwise, if $\rho(\beta)=1$, then  $\gamma=\overline{j}\cdot \alpha+\beta$ and 
$\tilde{\alpha}=\overline{j}\cdot \alpha=j^1\cdot \alpha.$  We find
$$F_{\gamma}^0(\tilde{\alpha})=(\tilde{\alpha}, \tilde{\alpha}+\gamma)=(\tilde{\alpha},\beta)=(\overline{j}\cdot \alpha, \beta)$$
while
$$F_{\gamma}^0(j\!:\!\tilde{\alpha})=(j^0\cdot\tilde{\alpha}, (j^1\cdot \tilde{\alpha})+\gamma)=(\alpha, \tilde{\alpha}+\gamma)=(\alpha, \beta).$$
So $F_{\gamma}^0$ carries the edge $j\#\tilde{\alpha}$ to $(\hat{j}\cdot \alpha, \beta).$

Now consider the alternate situation where $e=(\alpha, \hat{j}\cdot \beta)\in E(Q_{2q}),$ with $\alpha\in \mathbb{Z}_2^q, \beta\in B_q, j\in \{1, 2, \ldots , q\}.$  Then $\gamma=(j^0\cdot \alpha)+\beta\in B_q$.  We consider separately, subcases where $\rho(\alpha)=0$ versus where $\rho(\alpha)=1$.  If $\rho(\alpha)=0$, note that $j^0(\alpha)=\alpha$ and
$F_{\gamma}^0(\alpha)=(\alpha, \alpha+\gamma)=(\alpha, \beta)$ while 
$F_{\gamma}^0(j\!:\!\alpha)=(j^0\cdot \alpha, (j^1\cdot \alpha)+\gamma)
=(\alpha, \overline{j}\cdot \beta)$.  
So $F_{\gamma}^0$ carries edge $j\#\alpha$ to $(\alpha,\hat{j}\cdot \beta)$. 
If instead $\rho(\alpha)=1$, then $\alpha=j^1\cdot \alpha$ and 
$$F_{\gamma}^1(j\!:\!\alpha)=(j^1\cdot \alpha, (j^0\cdot \alpha)+\gamma)=(\alpha, \beta)$$
while
$$F_{\gamma}^1(\alpha)=(\alpha, \alpha+\gamma)=(\alpha, \overline{j}\cdot \beta)$$
so again, the edge $(\alpha, \hat{j}\cdot \beta)$ is covered.  This completes the proof of Proposition \ref{2*Qq divides Q2q}, and with it the proof of Theorem \ref{m*Qq divides Qmq}.  \hfill  $\Box$

%Section 4
\section{For $q$ odd, when does $P_m$ divide $Q_q$?}
\vspace{.2in}

  This section is devoted to answering the %above question (invert word order)
  question above.  We will show that for $q<2^{32}$, $P_m <_D Q_q$ if and only if $m\leq q$ and $m | q\cdot 2^{q-1}$. 
  
  We begin by citing from \cite{Ring1}, \cite{AS} and \cite{E-ZVE}, the fact that  $Q_{2n}$ has an edge-decomposition into Hamiltonian cycles.
  %Theorem 3
  \begin{theorem}\label{Ham decomp Q2n}  $Q_{2n}$ has an edge decomposition into $n$ copies of $C_{2^{2n}},$ {\it i.e.}, $C_{2^{2n}} <_D Q_{2n}$.  \hfill $\Box$
  \end{theorem}
  
  Because $P_{2^t} <_D C_{2^{2n}}$ when $t<2n$, an immediate corollary is  
  %Corollary 1
 \newtheorem{corollary}{Corollary}
    \begin{corollary}\label{P2^t divides Q2n}
  For any $t<2n, P_{2^t} <_D Q_{2n}.$ \hfill $\Box$
  \end{corollary}

We will also provide later a simple ``proof by formula" of the specialization of
Theorem \ref{Ham decomp Q2n} and Corollary \ref{P2^t divides Q2n}
to the case when $n$ is a power of 2.

Recall that the Cartesian product of two graphs $G$ and $G^{\,\prime}$, denoted $G\Box G^{\,\prime}$, is the graph whose vertex set is $V(G)\times V(G^{\,\prime})$ and whose edge set consists of pairs that are either $\{ (x_1, y), (x_2, y)\}$, where $\{x_1, x_2\}$ is an edge of $G$ and $y\in V(G^{\,\prime})$, or $\{ (x, y_1), (x, y_2)\}$ where $x \in V(G)$ and $\{y_1, y_2\}$ is an edge of $G^{\,\prime}$.
It should be clear that $Q_q\Box Q_{q^{\,\prime}}\simeq Q_{q+q^{\,\prime}}$.
To begin to relate Cartesian products and edge decompositions, we have

%Lemma 4
\begin{lemma}\label{Cartesian Products}

{\rm (a)} Let  $H$, $G$, $G^{\,\prime}$ be graphs.  If $H <_D G$ and $H <_D G^{\,\prime}$ then $H <_D G\Box G^{\,\prime}.$\\  
{\rm (b)} If $P_m <_D Q_q$ and  $P_m <_D Q_{q^{\,\prime}}$ then $P_m <_D Q_{q+q^{\,\prime}}$
\end{lemma}

\noindent {\em Proof.}  Part (a) is obvious because $E(G\Box G^{\,\prime})$ consists of $|\,V(G)\,|$ copies of $E(G^{\,\prime})$
 and $|\,V(G^{\,\prime})\,|$ copies of $E(G)$.  Part (b) is a specialization making use of $Q_q\Box Q_{q^{\,\prime}}\simeq Q_{q+q^{\,\prime}}$.  \hfill $\Box$

 %Prop 3
 \begin{proposition}\label{ind step P2^t divides Qq}
 Let $t<2n$ and suppose it can be shown that $P_{2^t} <_D Q_q$ for all odd integers $q$ in the range $2^t+1$ to $2^t+2n-1$.  Then $P_{2^t} <_D Q_q$ for all $q\geq 2^t$ ($q$ odd or even).
 \end{proposition}
 
 \noindent {\em Proof.}  If $q$ is even and $2^t\leq q$, put $n=q/2$.  It follows from Corollary \ref{P2^t divides Q2n} that $P_{2^t} <_D Q_q$ since $t<2^t\leq %
 q=2n$.  Now suppose $s$ is odd and $s>2n$, and the proposition holds for $q^{\,\prime}<2^t+s$, {\it i.e.}
$P_{2^t} <_D Q_{q^{\,\prime}}$ for $2^t\leq q^{\,\prime}<2^t+s.$    
Then $P_{2^t} <_D Q_{2^t+s-2n}$, and Lemma \ref{Cartesian Products}(b) and
Corollary \ref{P2^t divides Q2n} give us $P_{2^t} <_D Q_{2^t+s}$, {\it i.e.} the proposition holds for $q=2^t+s$ as well.  By induction on $s$ it holds for all odd $q>2^t.$
\hfill $\Box$

\vspace{.2in}

Proposition \ref{ind step P2^t divides Qq} shows that only a finite number of $Q_q$'s need to have path decompositions constructed, to infer that $P_{2^t} <_D Q_q$ for all $q\geq 2^t$.

The key idea for the construction is to extend paths shorter than length $2^t$ to paths of length $2^t$.  The object we utilize for doing this is defined next.

\begin{definition}\label{Definition DVOP}
Let $G$ be a graph.  A \underline{disjoint collection of vertex-originating}
\underline{paths of length k}, henceforth DVOP$[k]$, is a collection of paths of length $k$ $\{p_v: P_k \longrightarrow G\}$
indexed by $V(G)$, satisfying\\
{\rm (a)}  disjointness, {\it i.e.} $p_v(\hat{j})=p_{v^{\,\prime}}(\hat{j^{\,\prime}}) \Longrightarrow v=v^{\,\prime}\,{\rm and}\,j=j^{\,\prime}$, and\\
{\rm (b)}  $p_v(0)=v,$ {\em i.e.} each vertex originates one path.\\
Here, as above, the vertices of $P_k$ are taken to be $\{0, 1, \ldots , k\}$ and $\hat{j}$ denotes the edge joining $j-1$ and $j$.  
Clearly there are $k|\,V(G)\,|$ edges in the combined images of all the paths in a DVOP$[k]$.
\end{definition}

%Prop 4
\begin{proposition}\label{DVOP}
For $0\leq k\leq 3$ and $n\geq k$, there is an edge decomposition of $Q_{2n}$ into $n-k$ Hamiltonian cycles and a DVOP$[k]$.
\end{proposition}

\noindent {\em Proof.}  The case $k=0$ merely reiterates Theorem \ref{Ham decomp Q2n}.
For $k>0$ start with Theorem \ref{Ham decomp Q2n} giving an edge decomposition of
$Q_{2n}$ into $n$ copies of $C_{2^{2n}}$.  Call three of the cycles $C^{(1)}, C^{(2)}, C^{(3)}$
(or stop at $C^{(k)}$ if $k<3$), and choose a direction or orientation on each cycle.  
Define set bijections $h_i: V(Q_{2n}) \longrightarrow V(Q_{2n}), \, 1 \leq i \leq k$, by letting
$h_i(v)$ be the vertex reached by traveling one edge along $C^{(i)}$ from $v$,
in the chosen direction.  Let $p_v:P_k\rightarrow Q_{2n}$ be the path defined by:
$p_v(0)=v, p_v(1)=h_1(v), p_v(2)=h_2h_1(v), p_v(3)=h_3h_2h_1(v)$ (stop sooner if $k<3$).
This is clearly a graph map because $p_v(i-1)$ and $p_v(i)$ are connected by an edge in
$C^{(i)}$, and it obviously originates on $v$.  It is a path ({\it i.e.} an embedded copy of $P_k$) if the $k+1$ vertex images are all distinct.  Because adjacent vertices have opposite parity on a hypercube, $p_v(0) \neq p_v(1) \neq p_v(2) \neq p_v(3) \neq p_v(0)$.
To see that $p_v(0) \neq p_v(2)$ note that these two vertices are connected by distinct edges to $p_v(1)$ so must be distinct in $Q_{2n}$.  Likewise and $p_v(1) \neq p_v(3)$ because they
are connected by distinct edges to $p_v(2)$.   Finally, if two edges coincide, say
$p_v(\hat{i}) = p_w(\hat{j})$, then  because $p_v(\hat{i}) \in C^{(i)}$ while
$p_w(\hat{j}) \in C^{(j)}$ we have $C^{(i)} \cap C^{(j)} \neq \phi$ forcing $i=j$.
The fact that $p_v$ and $p_w$ follow the same chosen orientations of the $C^{(j)}$'s means
that $p_v(i)=p_w(i)$, and then bijectivity of the $h_j$'s leads to $v=w$.   \hfill  $\Box$
\vspace{.1in}

A key ingredient in the proof of the main theorem of this paper is 

\begin{definition}
\label{definition sharp}
Let $G$ be a graph and $m\geq 1.$  Let 
$$m\#G$$
denote the graph obtained by drawing two copies of $G$, $($call them $G^{\,\prime}$ and $G^{\,\prime\prime})$, and connecting each vertex $v^{\,\prime}\in V(G^{\,\prime})$ to the corresponding vertex $v^{\,\prime\prime}\in V(G^{\,\prime\prime})$ by a path of length $m$.
\end{definition}

%Lemma 5
\begin{lemma}
\label{props of sharp}
{\rm (a)}  $1\#G\simeq P_1\Box G.$\\
{\rm (b)}  $|\,V(m\#G)\,|=(m+1)|\,V(G)\,|.$\\
{\rm (c)}  $|\,E(m\#G)\,|=2|\,E(G)\,|+m|\,V(G)\,|.$
\end{lemma}

\noindent {\em Proof.}  Trivial.  \hfill  $\Box$

%Lemma 6
\begin{lemma}\label{If m is odd}
If $m$ is odd, $m\#Q_q <_D Q_{m+q}.$
\end{lemma}

\noindent {\em Proof.}  Denote a vertex of $Q_{m+q}$ as $(\alpha, \beta)$, where $\alpha\in V(Q_m)$ and $\beta\in V(Q_q)$.  Utilize the edge decomposition of $Q_m$ into copies of  $P_m$:
$$\{f_{\gamma}: P_m\longrightarrow Q_m |\,\gamma \in B_m\}$$
defined earlier.
Denote the embedded copies of $Q_q$ in $m\#Q_q$ as $Q_q^{\,\prime}$ and
$Q_q^{\,\prime\prime}$.  Denote the $j^{\textrm{th}}$ point on the edge of $m\#Q_q$ joining $\beta^{\,\prime}$ to $\beta^{\,\prime\prime}$ as $\beta_{\langle j\rangle}.$  Thus $\beta_{\langle 0\rangle}=\beta^{\,\prime}$ and
$\beta_{\langle m\rangle}=\beta^{\,\prime\prime}$.  Define for $\gamma\in B_m$
$$F_{\gamma}:  m\#Q_q\longrightarrow Q_{m+q}$$
by $F_{\gamma}(\beta_{\langle j\rangle})=(f_{\gamma}(j), \beta)$.  This is a collection of $2^{m-1}$ embeddings of $m\#Q_q$ into $Q_{m+q}$.  Since $|\,E(m\#Q_q)\,|=(m+q)2^q$, there are in all $2^{m-1}(m+q)2^q=(m+q)2^{m+q-1}$ edge images of all the $\{F_{\gamma}\}$.  Since $|\,E(Q_{m+q}\,|=(m+q)2^{m+q-1}$, the $\{F_{\gamma}\}$ are a decomposition into disjoint copies if their collective images are onto $E(Q_{m+q})$.  But this is easy, because an edge of $Q_{m+q}$ is either $(\hat{j}\cdot \alpha, \beta)$, where $\hat{j}\cdot \alpha$ is an edge of $Q_m$, or $(\alpha, \hat{k}\cdot \beta)$, where $\hat{k}\cdot \beta$ is an edge of $Q_q$.  Clearly $(\hat{j}\cdot \alpha, \beta)\in im (F_{\gamma})$ if $\hat{j}\cdot \alpha\in im (f_{\gamma})$.  The edge $(\alpha, \hat{k}\cdot \beta)$ equals $F_{\alpha}(\hat{k}\cdot\beta^{\,\prime})$ if $\alpha$ has even parity, and it equals
$F_{\gamma}(\hat{k}\cdot \beta^{\,\prime\prime})$ for 
$\gamma=\alpha+1^m$ if $\alpha$
has odd parity (note that $m$ must be odd for this to work).   \hfill  $\Box$

%Lemma 7
\begin{lemma}\label{DVOP and comp edge set}
Suppose $G$ has an edge decomposition into a DVOP$[k^{\,\prime}]$ and a complementary
edge set $E^{\,\prime}$, as well as an edge decomposition into a DVOP$[k^{\,\prime\prime}]$
and a complementary edge set $E^{\,\prime\prime}$. Then $m\#G$ has an edge decomposition
into $|\,V(G)\,|$ copies of $P_{k^{\,\prime}+m+k^{\,\prime\prime}}$ and one copy each of
$E^{\,\prime}$ and $E^{\,\prime\prime}$.
\end{lemma}

\noindent {\em Proof.}  Let $\{p_{v^{\,\prime}}:  P_{k^{\,\prime}}\longrightarrow G\}$ and
$\{p_{v^{\,\prime\prime}}:  P_{k^{\,\prime\prime}}\longrightarrow G\}$ be the DVOP's.  Simply
concatenate the paths $p_{v^{\,\prime}}\in G^{\,\prime}$, the path from $v^{\,\prime}$ to
$v^{\,\prime\prime}$ in $m\#G$, and the path $p_{v^{\,\prime\prime}}$ in
$G^{\,\prime\prime}$, to make the path
$\tilde{p_v}: P_{k^{\,\prime}+m+k^{\,\prime\prime}} \longrightarrow m\#G $.
Then $\{E(\tilde{p_v})\,|\,v\in V(G)\,\}\cup E^{\,\prime}\cup E^{\,\prime\prime}$ is an edge decomposition of $m\#G$.  \hfill  $\Box$

%Prop 5
\begin{proposition}\label{prop P4 and P8}

{\rm (a)}  $P_4 <_D Q_5.$\\
{\rm (b)}  $P_4 <_D Q_7.$\\
{\rm (c)}  $P_8 <_D Q_9.$\\
{\rm (d)}  $P_8 <_D Q_{11}.$
\end{proposition}

\noindent {\em Proof.}  (a).  Viewing $Q_5$ as $1\#Q_4$, apply Proposition \ref{DVOP} to obtain $E(Q_4)$ as a DVOP$[2]$ (with empty complementary set) and also write $E(Q_4)$ as a DVOP$[1]$ with a single $C_{16}$ as the complementary set.  We have an application of Lemma \ref{DVOP and comp edge set} with $k^{\,\prime}=2, k^{\,\prime\prime}=1, m=1.$  Thus $E(Q_5)$ decomposes into 16 copies of $P_4$ and one copy of $C_{16}$.  Since $C_{16}$ is 4 copies of $P_4$, we have shown that $P_4 <_D Q_5$.\\
(b)  Apply Proposition \ref{DVOP} for a DVOP$[1]$ and complementary set $C_{16}$ in $Q_4$, as well as an empty DVOP$[0]$ and complementary set consisting of 2 copies of $C_{16}$.  Then $3\#Q_4$ has an edge decomposition into 16 copies of $P_{1+3+0}=P_4$ and 3 copies of $C_{16}$, {\em i.e.} $P_4<_D 3\#Q_4$.  By Lemma \ref{If m is odd}, $P_4 <_D Q_7$.\\
(c)  The proof that $P_8<_D5\#Q_4$ likewise applies Lemma \ref{DVOP and comp edge set} with $k^{\,\prime}=2$ and $k^{\,\prime\prime}=1$, but now with $m=5$ so that paths have length $k^{\,\prime}+m+k^{\,\prime\prime}=8.$  It follows that $5\#Q_4$ has an edge decomposition into 16 copies of $P_8$ and one copy of $C_{16}$, hence a decomposition into 18 copies of $P_8$.  We have $P_8 <_D 5\#Q_4<_DQ_9$, so $P_8<_DQ_9$.\\
(d)  Similarly, $P_8 <_D 7\#Q_4$ (Lemma \ref{DVOP and comp edge set} with
$k^{\,\prime}=1, k^{\,\prime\prime}=0, m=7$),
and $7\#Q_4 <_D Q_{11}$.  \hfill   $\Box$

%Cor 2
\begin{corollary}\label{cor P4 and P8}
$P_4 <_D Q_q$ for all $q\geq 4$ and $P_8 <_D Q_q$ for all $q\geq 8$.
\end{corollary}

\noindent {\em Proof.}  This follows immediately from Propositions
\ref{ind step P2^t divides Qq} and \ref{prop P4 and P8}.  \hfill  $\Box$
\vspace{.1in}

\noindent {\bf Note}:  The first half of Corollary \ref{cor P4 and P8} was proved by an {\it ad hoc} method in [11].

%Lemma 8
\begin{lemma}\label{Q2n has DVOP[n]}
$Q_{2n}$ has a DVOP$[n] ($with empty complementary set$)$.
\end{lemma}

\noindent {\em Proof.}  Let $f_{\gamma}:  P_n\longrightarrow Q_n$ be defined as before, but without the restriction that $\gamma\in B_n$.  Write a vertex of $Q_{2n}$ as $(\alpha, \beta)$, where $\alpha\in V(Q_n), \beta\in V(Q_n)$.  Let 
$$p_{(\alpha, \beta)}(j)=\left \{ \begin{array}{ll} (f_{\alpha}(j), \beta) &{\rm  if}\, (\alpha, \beta)\, {\rm has}\,{\rm even}\,{\rm parity}\\
          							(\alpha, f_{\beta}(j)) & {\rm if}\, (\alpha, \beta) \,{\rm has}\,{\rm odd}\,{\rm parity}.
							\end{array} \right. $$
							
Then $p_{(\alpha, \beta)}:  P_n\longrightarrow Q_{2n}$ is a path and $p_{(\alpha, \beta)}(0)=(\alpha, \beta)$.  To prove disjointness, suppose two edges coincide, {\em e.g.}	$(\hat{j}\cdot \alpha, \beta)=(\hat{k}\cdot \alpha^{\,\prime}, \beta^{\,\prime}).$  Then obviously $\beta=\beta^{\,\prime}$ and $j=k$, and $\alpha$ and $\alpha^{\,\prime}$ differ at most in their $j^{th}$ coordinate.  But $\alpha+\beta$ and $\alpha^{\,\prime}+\beta^{\,\prime}=\alpha^{\,\prime}+\beta$ have the same parity, so $\alpha=\alpha^{\,\prime}$.  Likewise, if the coordinate that varies in the edge is among the last $n$ coordinates.  \hfill $\Box$	
\vspace{.1in}

To get results ilike Corollary \ref{cor P4 and P8} for $P_{16}$, we need to make use of $Q_8$.

%Lemma 9
\begin{lemma}\label{For 0 leq k leq 4, Q8}
For $0\leq k\leq 4, Q_8$ has a DVOP$[k]$ and a complementary set that consists of $4-k$ copies of $C_{256}$.
\end{lemma}					

\noindent {\em Proof.}  Use Lemma \ref{Q2n has DVOP[n]} (for $k=4$) and Proposition \ref{DVOP} (for $k<4$).  \hfill  $\Box$

%Lemma 10
\begin{lemma}\label{For t=4, 5, 6, 7, and for s=1, 3, 5, 7}
For $t=4, 5, 6, 7, $ and for $s=1, 3, 5, 7,$ 
$P_{2^t} <_D Q_{2^t+s}.$
\end{lemma}

\noindent {\em Proof.}  Let $m=2^t+s-8.$  Put $k^{\,\prime}=0$ if $s>4$ and put $k^{\,\prime}=4$ if $s<4$.  Put $k^{\,\prime\prime}=8-s-k^{\,\prime}$, which will equal either $1$ or $3$ depending on $s$.
Apply Lemmas \ref{For 0 leq k leq 4, Q8}, \ref{DVOP and comp edge set}, and
\ref{If m is odd} to deduce that $P_{2^t} <_D m\#Q_8<_D Q_{m+8}=Q_{2^t+s}$.  \hfill $\Box$

%\begin{flushright} $\Box$  \end{flushright}

%Corollary 3
\begin{corollary}\label{P2^t divides Qq if and only if}
For $t=4, 5, 6, 7,$  
$P_{2^t} <_D Q_q\, {\rm if}\,\,{\rm and}\,{\rm only}\,{\rm if}\,q\geq 2^t.$
\end{corollary} 

\noindent {\em Proof.}  This follows from Proposition \ref{ind step P2^t divides Qq} (put $n=4$)
and Lemma \ref{For t=4, 5, 6, 7, and for s=1, 3, 5, 7}.  \hfill $\Box$
\vspace{.1in}

The natural generalization of Proposition \ref{prop P4 and P8} and Lemma \ref{For t=4, 5, 6, 7, and for s=1, 3, 5, 7} is 
%Prop 6
\begin{proposition}\label{Suppose for each odd k between 0 and 2^{r-1}}
Let $r\geq 2$.  Suppose that for each odd number $k$ between $0$ and $2^{r-1}, Q_{2^r}$ has an edge decomposition into a DVOP$[k]$ and $2^{r-1}-k$ copies of $C_{2^{2^r}}$.  Then for any $t$ in the range $2^{r-1}\leq t<2^r,$\\
{\rm (a)}  $P_{2^t} <_D Q_{2^t+s}$ for $s$ odd, $1\leq s\leq 2^r-1;$ and\\
{\rm (b)}  $P_{2^t} <_D Q_q$ if and only if $q\geq 2^t.$
\end{proposition}

\noindent {\em Proof.}  (b) follows from (a) and Proposition \ref{ind step P2^t divides Qq} .\\
The proof of (a) is a straightforward generalization of the proof Lemma  \ref{For t=4, 5, 6, 7, and for s=1, 3, 5, 7}.  Put $m=2^t+s-2^r$.  Let $k^{\,\prime}$ denote either 0 if $s>2^{r-1}$, or $k^{\,\prime}=2^{r-1}$ if $s<2^{r-1}$.  Put $k^{\,\prime\prime}=2^r-s-k^{\,\prime}$.  Then $k^{\,\prime\prime}$ is an odd number between 0 and $2^{r-1}$.  Note that
$k^{\,\prime}+m+k^{\,\prime\prime}=2^t.$  Apply Lemmas  \ref{If m is odd}, \ref{DVOP and comp edge set} and \ref{Q2n has DVOP[n]} to obtain that
$P_{2^t} <_D m\#Q_{2^r} <_D Q_{m+2^r}=Q_{2^t+s}.$         \hfill  $\Box$
\vspace{.1in}

It is tempting to look for DVOP$[k]$'s in $Q_{2^r}$ by stitching together paths consisting of one edge from each of $k$ Hamiltonian cycles.  The resulting edge sets are vertex-originating and are disjoint.  What seems to be hard, is to prove they embed paths, {\it i.e.} that no vertex is repeated.  Proposition \ref{DVOP} exploited the fact that any map from $P_3$ to $Q_n$ having distinct edge images embeds a path, because there are no loops of length 2 or 3 in $Q_n$.  Generalizing Proposition \ref{DVOP} takes some work, and we have succeeded only for $r=4$ and $5$.
\vspace{.1in}

\noindent {\bf Construction}
For $r\geq 1$,  we define a set of $2^{r-1}$ cycles in $Q_{2^r}$ indexed by
$\delta=(d_1, \ldots , d_{r-1})\in \mathbb{Z}_2^{r-1}$, denoted
$g_{\delta}:  C_{2^{2^r}}\longrightarrow Q_{2^r}$, as follows.
For $r=1$ there is just one cycle, denoted $g:  C_4\longrightarrow Q_2$, which traces
the unique cycle starting at $0^2$, {\em i.e.} $g(0)=(0, 0); g(1)=(1, 0); g(2)=(1, 1); g(3)=(0, 1)$.
The vertices of $C_n$ are identified with the integer range $[0,n-1]$ viewed modulo $n$.
For $r=2$ define $g_0, g_1:  C_{16}\longrightarrow Q_4$ by
$$g_0(4u+v)=(g(v-u), g(u));  g_1(4u+v)=(g(u), g(v-u)); \,{\rm where}\,\, 0\leq u, v\leq 3.$$
Then $\{g_0, g_1\}$ is a set of two maps from $[0, 15]$ to $Q_4$ indexed by $\mathbb{Z}_2$, and their images turn out to be edge-disjoint cycles that partition $E(Q_4)$.
Now let $r \geq 2$ and suppose the
$\{g_{\delta}:  C_{2^{2^r}}\longrightarrow Q_{2^r}\,|\,\delta \in \mathbb{Z}_2^{r-1} \}$
have been defined.  The vertices of $Q_{2^{r+1}}$ will be identified with $Q_{2^r}\Box Q_{2^r}$ and may be written as $(\alpha,\beta)$, where $\alpha, \beta\in V(Q_{2^r})$.
For $\delta=(d_1, \ldots , d_{r-1})\in \mathbb{Z}_2^{r-1}$, let
$\delta 0$ (respectively $\delta 1$) denote
$(d_1, \ldots ,d_{r-1}, 0)$ (respectively $(d_1, \ldots , d_{r-1}, 1)$) $\in \mathbb{Z}_2^{r}$.
Define the cycles  $\{g_{\delta 0} : C_{2^{2^{r+1}}}\longrightarrow Q_{2^{r+1}}\}$ and
$\{g_{\delta 1} : C_{2^{2^{r+1}}}\longrightarrow Q_{2^{r+1}}\}$ by these formulas: 
$$g_{\delta 0}(2^{2^r}u+v)=\big(g_{\delta} (v-u), g_{\delta}(u)\big),$$
$$g_{\delta 1}(2^{2^r}u+v)=(g_{\delta} (u), g_{\delta}(v-u)),$$
for $u, v\in [0, 2^{2^r}-1]$.
Taken together, $\{g_{\delta 0}\} \cup \{g_{\delta 1}\}$ is a set of $2^r$ cycles
in $Q_{2^{r+1}}$, indexed by $\mathbb{Z}_2^{r}$, that comprises the construction for $r+1$.

%Prop 7
\begin{proposition}\label{Ham decomp E(Q2^r)}  For
$r\geq 1, \{g_{\delta}: C_{2^{2^r}}\longrightarrow  Q_{2^r} \,|\,\delta\in \mathbb{Z}_2^{r-1})\}$ is a decomposition of $E(Q_{2^r})$ into a disjoint collection of $2^{r-1}$ Hamiltonian cycles.
\end{proposition}

\noindent {\em Proof.}  It is true for $r=1$, where the collection consists of the singleton
$\{g: C_4\longrightarrow Q_2\}$.  Assuming it is true for $r$, we prove it for $r+1$.
First note that $g_{\delta 0}: [0, 2^{2^{r+1}}-1]\longrightarrow Q_{2^{r+1}}$ is bijective because
$g_{\delta}$ is, and likewise for $g_{\delta 1}$.  Second, $g_{\delta 0}$ is indeed a cycle, {\it i.e.}
it carries edges of $C_{2^{2^{r+1}}}$ to a closed chain of edges of $Q_{2^{r+1}}$.
If $v\neq 2^{2^r}-1, g_{\delta 0}$ maps the edge from $2^{2^r}u+v$ to $2^{2^r}u+v+1$ into
the edge from $(g_{\delta}(v-u), \beta)$ to $(g_{\delta}(v-u+1), \beta)$ in $Q_{2^{r+1}}$, where
$\beta = g_{\delta}(u).$   If $v=2^{2^r}-1,$ the edge from $2^{2^r}u+v$ to $2^{2^r}(u+1)$
[if $u=2^{2^r}-1$, we mean the edge of $C_{2^{2^{r+1}}}$ from $2^{2^{r+1}}-1$ to $0$]
goes to the edge from $(\alpha, g_{\delta}(u))$ to $(\alpha, g_{\delta}(u+1)$, where
$\alpha=g_{\delta}(v-u)$.  The proof for $g_{\delta 1}$ is identical.  Finally, the cycles are
disjoint and cover $E(Q_{2^{r+1}})$.  To see this, since the total number of edge images
of all the $\{g_{\delta 0}\, {\rm and}\,g_{\delta 1}\}$ is $(2^r)(2^{2^{r+1}})$, which equals
$|\,E(Q_{2^{r+1}})\,|,$ it suffices to show that every edge of $E(Q_{2^{r+1}})$ is covered.
Suppose an edge of $E(Q_{2^{r+1}})$ is $(\hat{j}\cdot \alpha, \beta)$; the proof for
$(\alpha, \hat{j}\cdot \beta)$ is equivalent.  Recursively, we know that
$\hat{j}\cdot \alpha$ is the image under some $g_{\delta},$ of the edge of
$C_{2^{2^r}}$ from $w$ to $w+1$ for some $w\in [0, 2^{2^r}-1].$ 
Let $u=g_{\delta}^{-1}(\beta)$ and let $v=w+u$ if $w+u<2^{2^r}$, and $v=w+u-2^{2^r}$
otherwise. Then $u, v\in [0, 2^{2^r}-1].$ If $v\neq 2^{2^r}-1$ then $(\hat{j}\cdot\alpha,\beta)$ is the image under $g_{\delta 0}$ of the edge from $2^{2^r}u+v$ to $2^{2^r}u+v+1.$
If $v=2^{2^r}-1$, then $(\hat{j}\cdot \alpha, \beta)$ is the image under $g_{\delta 1}$
of the edge from $2^{2^r}u+v$ to $2^{2^r}(u+1).$  \hfill  $\Box$
\vspace{.1in}

\noindent {\bf Note}:  Proposition \ref{Ham decomp E(Q2^r)} gives the same type of cycle decomposition as Theorem \ref{Ham decomp Q2n}.  However, we want the notations and formulas for $g_{\delta}$ for other purposes, so we have offered this alternative proof for the situation
where $2n=2^r$.  This article only makes use of cycle decompositions for hypercubes that are
$Q_{2^t}$, so having Proposition \ref{Ham decomp E(Q2^r)} also keeps the the principal results entirely self-contained.
\vspace{.1in}

\begin{definition}\label{Definition 1-value}
Let $n\geq 1$  Define the \underline{$1$-value} of a vertex $\alpha=(a_1, \ldots, a_{2n})\in V(Q_{2n})$ by 
$$\rho_1(\alpha)=\sum_{i=1}^n g^{-1}(a_{2i-1}, a_{2i})  \pmod 4,$$
where $g:  [0, 3]\longrightarrow Q_2$ is the map defined above.\\
Likewise the \underline{$2$-value} of a vertex $\alpha=(a_1,\ldots, a_{4n})\in Q_{4n}$ is\\
$$\rho_2(\alpha)=\sum_{i=1}^n g_0^{-1}(a_{4i-3}, a_{4i-2}, a_{4i-1}, a_{4i})  \pmod {16},$$
where $g_0: [0, \ldots , 15] \longrightarrow Q_4$ is the map defined above.
\end{definition}

%Lemma 11
\begin{lemma}\label{list of 1-values}
Let $r\geq 1$ and let $\delta\in \mathbb{Z}_2^{r-1}.$
Then $\rho_1(g_{\delta}(w))\equiv w \pmod {4}.$
That is, as any of the $2^{r-1}$ Hamiltonian cycles of $Q_{2^r}$ defined above is traversed, the
$1$-values of the vertices visited increases by $+1 \pmod {4}$ at each step.  The list of
$1$-values seen along any of these Hamiltonian cycles is:  $0, 1, 2, 3, 0, 1, 2, 3, 0, \ldots$
\end{lemma}

\noindent {\em Proof.}  It is true for the unique cycle of $Q_2$, {\it i.e.} for $r=1$.
Assuming it is true for $r$, let $\delta\in \mathbb{Z}_2^{r-1}$ and let $u, v\in [0, 2^{2^r}-1]$.
Then in $Q_{2^{r+1}}$,
$$\rho_1(g_{\delta 0}(2^{2^r}u+v))=\rho_1(g_{\delta}(v-u), g_{\delta}(u))$$
$$=\rho_1g_{\delta}(v-u)+\rho_1g_{\delta}(u)\equiv (v-u)+u\equiv 2^{2^r}u+v \pmod{4}.$$
Likewise for $\rho_1(g_{\delta 1} (2^{2^r}u+v))$, confirming that the formula holds for $r+1$.
 \hfill  $\Box$

%Prop 8
\begin{proposition}\label{Let r geq 4 and 0 leq k leq 7}
Let $r\geq 4$ and $0\leq k\leq 7$.  There is an edge decomposition of $Q_{2^r}$ into
$2^{r-1}-k$ Hamiltonian cycles $($copies of $C_{2^{2^r}})$ and a DVOP$[k]$.
\end{proposition}

\noindent {\em Proof.}  For $\delta\in \mathbb{Z}_2^{r-1}$, let
$h_{\delta}: Q_{2^r}\longrightarrow Q_{2^r}$ be defined by
$h_{\delta}(\alpha)=g_{\delta}(g_{\delta}^{-1}(\alpha)+1)$, {\it i.e.}, $h_{\delta}$ consists of
advancing one edge along the Hamiltonian cycle indexed by $\delta$.
Because $r \geq 4$ there are at least 8 distinct $h_{\delta}$'s, and we abuse notation slightly
to let $h_1, \ldots, h_7$ denote any 7 distinct $h_{\delta}$'s from that set.

For each $\alpha \in V(Q_{2^r})$ define a path $p_{\alpha}: P_k\longrightarrow Q_{2^r}$ by:
$p_{\alpha}(0)=\alpha$; $p_{\alpha}(1)=h_1(\alpha)$;
$p_{\alpha}(2)=h_2^{-1}(p_{\alpha}(1))$; $p_{\alpha}(3)=h_3(p_{\alpha}(2))$;
$p_{\alpha}(4)=h_4(p_{\alpha}(3))$; $p_{\alpha}(5)=h_5(p_{\alpha}(4))$;
$p_{\alpha}(6)=h_6^{-1}(p_{\alpha}(5))$; $p_{\alpha}(7)=h_7(p_{\alpha}(6))$
(stop sooner if $k<7$).  Then all the edges used are distinct because the $h_{\delta}$'s are
bijections and the Hamiltonian cycles are disjoint.  Obviously $p_{\alpha}$ originates on
$\alpha$.  To see that the set of (up to) $8$ vertices visited on this $p_{\alpha}$ are distinct,
look at their list of 1-values.  Letting $\nu=\rho_1(\alpha)$, the list is $\pmod 4$:
$$\nu, \nu+1, \nu, \nu+1, \nu+2, \nu+3, \nu+2, \nu+3.$$
The list derives from Lemma \ref{list of 1-values} and the fact that motion along the second and
sixth cycles goes opposite to its usual orientation (by using $h_2^{-1}$ and $h_6^{-1}$).  The
list proves that the image of $p_{\alpha}$
is a path since a repeated vertex would repeat its $1$-value, but repeats of $1$-values
are always two vertices apart and there are no $2$-cycles on a hypercube.   \hfill  $\Box$ 

%Cor 4
\begin{corollary}\label{For 8 leq t leq 15}
For $8\leq t\leq 15, P_{2^t} <_D Q_q$ for $q\geq 2^t.$
\end{corollary}

\noindent {\em Proof.}  Propositions \ref{Let r geq 4 and 0 leq k leq 7} shows that the
hypothesis of Proposition \ref{Suppose for each odd k between 0 and 2^{r-1}}
holds for $r=4$, and thus so does the conclusion.
\hfill    $\Box$

%Lemma 12
\begin{lemma}\label{let g_0, g_1 : [0, 15] longrightarrow Q4}
Let $g_0, g_1 : [0, 15] \longrightarrow Q_4$ be the cycles defined above, and let the subscript
$0\delta$ mean $(0, d_1, \ldots , d_{r-1})$ for $\delta=(d_1, \ldots ,, d_{r-1})$ and likewise for
$1\delta$.  Similarly define $0\delta 0, 0\delta 1, 1\delta 0, \,{\rm and}\,1\delta 1.$\\
{\rm (a)} For $w \in [0,15], \,\rho_2(g_0(w))\equiv w \pmod {16}\,$ and 
$\,\rho_2(g_1(w))\equiv 5w \pmod {8}$.\\
{\rm (b)}  For any $r\geq 2,$ and any $\delta\in \mathbb{Z}_2^{r-2}$,
$$\begin{array}{lcrll}
\rho_2(g_{0\delta}(w))&\equiv &w &\pmod {16}, &\,{\rm and}\\
\rho_2(g_{1\delta}(w))&\equiv &5w &\pmod 8.& \end{array}$$
{\rm (c)}  For $r \geq 2$, if $2$-values of vertices are reduced modulo $8$, then advancing
one step along a cycle indexed by any $0\delta$ adds $+1 \pmod 8$ while advancing
one step along any cycle indexed by a $1\delta$ adds $+5 \,\, ({\rm or} -3)\pmod 8$.
\end{lemma}

\noindent {\em Proof.}  Part (a) is easiest proved simply by writing down the values of
$g_0(w)$ and $g_1(w)$ for $w=0, \ldots , 15.$
Part (c) follows easily from part (b).
Part (a) shows that (b) is true for $r=2$.  If (b) holds at $r$, then at $r+1$ we have
$\rho_2(g_{0\delta 0}(2^{2^r}u+v))=\rho_2(g_{0\delta}(v-u), g_{0\delta}(u))
=\rho_2(g_{0\delta}(v-u))+\rho_2(g_{0\delta}(u))\equiv (v-u)+u=v
\equiv 2^{2^r}u+v \pmod {16}$.
The proofs for $\rho_2g_{0\delta 1}, \rho_2g_{1\delta 0},  {\rm and}\, \rho_2g_{1\delta 1}$
are similar.  \hfill  $\Box$ 

%Prop 9
\begin{proposition}\label{r geq 5 and 0 leq k leq 15}
Let $r\geq 5$ and $0\leq k\leq 15$.  There is an edge decomposition of $Q_{2^r}$ into
$2^{r-1}-k$ Hamiltonian cycles $($copies of $C_{2^{2^r}})$ and a DVOP$[k]$.
\end{proposition}

\noindent {\em Proof.}  As in the proof of Proposition \ref{Let r geq 4 and 0 leq k leq 7},
define $h_{\delta}: Q_{2^r}\longrightarrow Q_{2^r}$ by
$h_{\delta}(\alpha)=g_{\delta}(g_{\delta}^{-1}(\alpha)+1).$  It is a vertex bijection with the
property that $\alpha$ and $h_{\delta}(\alpha)$ are joined by an edge that belongs to the
$\delta^{\textrm{th}}$ Hamiltonian cycle.  Let $h_1, \ldots ,h_8$ denote any eight of these
of the form $h_{0\delta}$, and let $\overline{h}_1, \ldots , \overline{h}_7$ denote any seven
of these that are of the form $h_{1\delta}$ (because $r \geq 5$ there are enough
$h_{0\delta}$'s and $h_{1\delta}$'s for these to be chosen to be distinct).
Starting with a vertex $\alpha \in V(Q_{2^r}),$ define a sequence of (up to) 15 more vertices
by successively applying this list of functions:
$$h_1, h_2^{-1}, h_3, h_4, \overline{h}_1, \overline{h}_2^{-1}, \overline{h}_3, \overline{h}_4, h_5, h_6^{-1}, h_7, h_8, \overline{h}_5, \overline{h}_6^{-1}, \overline{h}_7.$$ 
These functions alter $2$-values in a consistent way:  $h_i, h_i^{-1}, \overline{h}_i, {\rm and}\,\overline{h}_i^{-1}$ alter the $2$-value by $+1, -1, -3, {\rm and}\,+3 \pmod{8}$ respectively.
Putting $\nu = \rho_2(\alpha),$ the resulting path $p_{\alpha}$ visits sequentially vertices
with the following list of $2$-values $\pmod{8}$:
$$\nu, \nu+1, \nu, \nu+1, \nu+2, \nu+7, \nu+2, \nu+7, \nu+4, \nu+5, \nu+4, \nu+5, \nu+6, \nu+3, \nu+6, \nu+3.$$
Vertices sharing a $2$-value $\pmod{8}$ are always just two edges apart, again making a
repeated vertex impossible.  Thus $p_{\alpha}: P_k\longrightarrow Q_{2^r}$ is a path and we
have defined a DVOP$[k]$.  The complementary set consists of the $2^{r-1}-k$ unused
Hamiltonian cycles.  \hfill   $\Box$

%Cor 5
\begin{corollary}\label{16 leq t leq 31}
For $16\leq t\leq 31, P_{2^t} <_D Q_q$ for $q\geq 2^t.$
\end{corollary}

\noindent {\em Proof.}  By Proposition \ref{r geq 5 and 0 leq k leq 15}, the hypothesis of Proposition \ref{Suppose for each odd k between 0 and 2^{r-1}} holds for $r=5$.  \hfill  $\Box$
\vspace{.2in}

Finally, let us put these results together to see what we can say about the motivating question:  when does $P_m$ divides $Q_q$ for $q$ odd?

%Theorem 4
\begin{theorem}\label{necessary condition}
Let $q$ be odd.  A necessary condition for $P_m$ to divide $Q_q$ is that
$m\leq q$ and $m\,|\,q\cdot 2^{q-1}$.
\end{theorem}

\noindent {\em Proof.}  That $m\,|\,q\cdot 2^{q-1}$ is obvious since $|\,E(Q_q)\,|$ must be a multiple of $|\,E(P_m)\,|$.  Because every vertex of $Q_q$ has odd degree, at least one path must start or end there.  Each path provides just two ``starts" or ``ends", and there are $q\cdot 2^{q-1}/m$ paths, hence $2(q\cdot 2^{q-1}/m)\geq |\,V(Q_q)\,|=2^q.$  This reduces to $q\geq m.$  \hfill   $\Box$ 
\vspace{.2in}

\noindent {\bf Conjecture}.  The above necessary condition is also sufficient.
For $q$ odd, $P_m <_D Q_q$ if $m\leq q$ and $m\,|\,q\cdot 2^{q-1}.$

%Theorem 5`
\begin{theorem}\label{main theorem}
The conjecture is true for $q<2^{32}$.
\end{theorem}

\noindent {\em Proof.}  Let $d=gcd(m, q)$.  Because $q$ is odd, $d$ is odd.
Consider the cases $d=1$ and $d>1$ separately.
If $d=1, \,m\,|\,2^{q-1}$ so $m$ is a power of $2$.
Let $2^t$ be the largest power of $2$ that is smaller than $q$.
Since $m\leq q, \,m\,|\,2^t.$  So $P_m<_D P_{2^t}$ and we only have to show that $P_{2^t} <_D Q_q$.  For $q<2^{32}, t<32$ so this is true by Corollaries \ref{cor P4 and P8} $-$ \ref{16 leq t leq 31} (or by the trivial case $P_2  <_D Q_q$ for $q\geq 2$).

Now suppose $d>1\,$: we reduce the case of $d=1$.  Let $m^{\,\prime} = m/d$ and let
$q^{\,\prime} = q/d.$  Then $m^{\,\prime} \,|\, q^{\,\prime} \cdot 2^{q-1}.$
But $m^{\,\prime}$ and $q^{\,\prime}$
are relatively prime so $m^{\,\prime} \,|\, 2^{q-1}$, making $m^{\,\prime}$ a power of $2$.
Since $m^{\,\prime} \leq q^{\,\prime} \leq 2^{q^{\,\prime} -1}$, we see that
$m^{\,\prime} \,|\, 2^{q^{\,\prime} -1} \,|\, q^{\,\prime} \cdot 2^{q^{\,\prime} -1}.$
Then $m^{\,\prime}$ and $q^{\,\prime}$ are less than $2^{32}$ and are relatively prime
so fall under the previous case, hence $P_{m^{\,\prime}} <_D Q_{q^{\,\prime}},$ {\it i.e.}
$P_{m/d} <_D Q_{q/d}.$ Apply Theorem \ref{m*Qq divides Qmq} to see that 
$P_m=d*P_{m/d} <_D d*Q_{q/d} <_D Q_q.$
\hfill  $\Box$

 \end{document}